\documentclass[a4paper]{amsart}
\usepackage{multicol,ifthen}
\usepackage{amssymb} 

\begin{document}

\newcommand{\F}{\mathcal{F}}
\newcommand{\R}{\mathbb R}
\newcommand{\N}{\mathbb N}
\newcommand{\C}{\mathbb C}  
\newcommand{\h}[2]{\mbox{$ \widehat{H^{#1}_{#2}}$}}
\newcommand{\hh}[3]{\mbox{$ \widehat{H^{#1}_{#2, #3}}$}} 
\newcommand{\n}[2]{\mbox{$ \| #1\| _{ #2} $}} 
\newcommand{\x}{\mbox{$X^r_{s,b}$}} 
\newcommand{\xx}{\mbox{$X_{s,b}$}}
\newcommand{\X}[3]{\mbox{$X^{#1}_{#2,#3}$}} 
\newcommand{\XX}[2]{\mbox{$X_{#1,#2}$}}
\newcommand{\q}[2]{\mbox{$ {\| #1 \|}^2_{#2} $}}
\newcommand{\e}{\varepsilon}
\newtheorem{lemma}{Lemma} 
\newtheorem{kor}{Corollary} 
\newtheorem{satz}{Theorem}
\newtheorem{prop}{Proposition}

\title[mKdV equation in almost critical $\widehat{H^r_s}$-spaces]{Local well-posedness for the modified KdV equation in almost critical $\widehat{H^r_s}$-spaces}

\author[A. Gr\"unrock and L. Vega]{Axel Gr\"unrock \\Fachbereich C:  Mathematik/Naturwissenschaften\\Bergische
Universit\"at Wuppertal, D-42097 Wuppertal, Germany \\
Axel.Gruenrock@math.uni-wuppertal.de \\ \hfill \\
\and \\ \hfill \\
Luis Vega \\ Departamento de Matematicas \\ Universidad del Pais Vasco, 48080 Bilbao, Spain\\
luis.vega@ehu.es}

\date{}

\begin{abstract}
We study the Cauchy problem for the modified KdV equation
\[u_t + u_{xxx} + (u^3)_x = 0, \hspace{2cm} u(0)=u_0\]
for data $u_0$ in the space $\h{r}{s}$ defined by the norm
\[\n{u_0}{\h{r}{s}} := \n{\langle \xi \rangle ^s\widehat{u_0}}{L^{r'}_{\xi}}.\]
Local well-posedness of this problem is established in the parameter range $2 \ge r >1$,
$s \ge \frac{1}{2} - \frac{1}{2r}$, so the case $(s,r)=(0,1)$, which is critical in view of scaling
considerations, is almost reached. To show this result, we use an appropriate variant of the Fourier 
restriction norm method as well as bi- and trilinear estimates for solutions of the 
Airy equation.
\end{abstract}

\maketitle

\section{Introduction and main result}

In this paper we study the local well-posedness (LWP) of the Cauchy problem for the modified KdV equation
\begin{equation}\label{10}
u_t + u_{xxx} + (u^3)_x = 0, \hspace{2cm} u(0)=u_0, \hspace{2cm}x\in \R.
\end{equation}
As long as data $u_0$ in the classical Sobolev spaces $H^s_x$ are considered, this problem is known to be 
well-posed for $s \ge \frac{1}{4}$ and ill-posed (in the $C^0$ - uniform sense) for $s < \frac{1}{4}$. Both,
the positive and the negative result, were shown by Kenig, Ponce, and the second author, see \cite[Theorem 2.4]{KPV93}
and \cite[Theorem 1.3]{KPV01}, respectively. The situation remains the same, when the defocusing modified
KdV equation, i. e. (\ref{10}) with a negative sign in front of the nonlinearity, is considered. In this case 
the proof of the well-posedness result remains identically valid, while the ill-posedness result here is due 
to Christ, Colliander and Tao, cf. \cite[Theorem 4]{CCT03}. In both cases the standard scaling argument suggests 
LWP for $s > -\frac{1}{2}$, so - on the $H^s_x$-scale - there is a considerable gap of $\frac{3}{4}$ derivatives 
between the scaling prediction and the optimal LWP result.

This gap could be closed partially by the first author in \cite{G04}, where data in the spaces $\h{r}{s}$ 
are considered, which are defined by the norms
\[\n{u_0}{\h{r}{s}} := \n{\langle \xi \rangle ^s\widehat{u_0}}{L^{r'}_{\xi}},\]
where $\widehat{u_0}$ denotes the Fourier transform of $u_0$, $\langle \xi \rangle = (1+\xi^2)^{\frac{1}{2}}$ and $\frac{1}{r}+\frac{1}{r'}=1$. The choice of these 
norms was motivated by earlier work of Cazenave, Vilela and the second author on nonlinear Schr\"odinger equations, 
see \cite{CVV01}, yet another alternative class of data spaces has been considered in \cite{VV01}.

The main result in \cite{G04} was LWP for (\ref{10}) in the parameter range $2\ge r > \frac{4}{3}$, $s\ge s(r):=\frac{1}{2}-\frac{1}{2r}$,
which coincides for $r=2$ with the optimal result on the $H^s_x$-scale. The proof used an appropriate variant 
of Bourgain's Fourier restriction norm method, cf. \cite{B93}. Especially the function spaces $\x$, defined by
\[\n{f}{\x}:= \left(\int d \xi d \tau \langle \xi \rangle^{sr'}\langle \tau - \xi^3\rangle^{br'} |\hat{f}(\xi , \tau)|^{r'} \right) ^{\frac{1}{r'}},\,\,\,\frac{1}{r}+\frac{1}{r'}=1\]
were utilised, as well as the time restriction norm spaces
\[\x(\delta) := \{f = \tilde{f}|_{[-\delta,\delta] \times \R} : \tilde{f} \in \x\}\]
with norm
\[\n{f}{\x(\delta)}:= \inf \{ \n{\tilde{f}}{\x} : \tilde{f}|_{[-\delta,\delta] \times \R} =f\}  .\]
A key estimate in \cite{G04} was the following Airy-version of the Fefferman-Stein-estimate (cf. \cite{F70} and \cite[Corollary 3.6]{G04})
\begin{equation}\label{20}
\n{e^{-t\partial^3 u_0}}{L^{3r}_{xt}} \le c \n{I^{-\frac{1}{3r}}u_0}{\widehat{L^r_{x}}},\hspace{1cm}r>\frac{4}{3}.
\end{equation}
Here and below $I$ ($J$) denotes the Riesz (Bessel) potential operator of order $-1$ and $\widehat{L^r_{x}}=\h{r}{0}$. This 
estimate fails to be true for $r \le \frac{4}{3}$, which explains the restriction $r > \frac{4}{3}$ in \cite{G04}.

It is the aim of the present paper to show, how this difficulty can be overcome by using bi- and trilinear estimates 
for solutions of the Airy equation (instead of linear and bilinear ones). This allows us to extend the LWP result 
for (\ref{10}) to the parameter range $2\ge r > 1$, $s\ge s(r)$. More precisely, the following theorem is the main result 
of this paper.

\begin{satz}\label{t1}
Let $2 \ge r > 1$, $s \ge s(r)=\frac{1}{2} - \frac{1}{2r}$ and $u_0\in \h{r}{s}$. Then there exist $b>\frac{1}{r}$, 
$\delta=\delta(\n{u_0}{\h{r}{s}})>0$ and a unique solution $u \in \x (\delta)$ of (\ref{10}). This solution is 
persistent and the flow map $S:u_0 \mapsto u,\,\,\,\h{r}{s} \rightarrow \x(\delta_0)$ is locally Lipschitz 
continuous for any $\delta_0 \in (0,\delta)$.
\end{satz}

Theorem \ref{t1} is sharp in the sense that, for given $r \in (1,2]$, we have ill-posedness in the $C^0$-uniform 
sense for $\frac{1}{r}-1 <s <s(r)$. This can be seen by using the counterexample from \cite{KPV01}, as it was discussed 
in \cite[section 5]{G04}. Combined with scaling considerations - observe that $\h{r}{s}$ scales like $H_x^{\sigma}$,
if $s-\frac{1}{r}=\sigma-\frac{1}{2}$ - this shows, that the case $(s,r)=(0,1)$ becomes critical in our setting 
and that our result covers the whole subcritical range. Unfortunately, our argument breaks down - even for small data - 
in the critical case, and we must leave this as an open problem. Notice, however, that for specific data
\[u_0= a \,\,\delta + \mu \,\,p.v. \frac{1}{x} \hspace{2cm}(a,\mu \,\,\,\,\mbox{small})\]
of critical regularity the existence of global solutions of (\ref{10}) was shown in \cite[Theorem 1.2]{PV05}.
By the general LWP Theorem \cite[Theorem 2.3]{G04} the proof of the following estimate is sufficient to establish 
Theorem \ref{t1}.

\begin{satz}\label{t2} Let $2 \ge r > 1$ and $s \ge s(r)=\frac{1}{2} - \frac{1}{2r}$. Then for all $b' < 0$ and $b > \frac{1}{r}$ the estimate
\begin{equation}\label{40}
\n{\partial _x (\prod_{i=1}^3 u_i)}{\X{r}{s}{b'}}\le c \prod_{i=1}^3 \n{u_i}{\x}
\end{equation}
holds true.
\end{satz}

\underline{Remarks:}
\begin{itemize}
\item[i)] (On the lifespan of local solutions) Using \cite[Lemma 5.2]{G05}, we have for $u_1$, $u_2$, $u_3$ supported in
$[-\delta,\delta] \times \R$ ($0<\delta \le 1$) the estimate
\[\n{\partial _x (\prod_{i=1}^3 u_i)}{\X{r}{s}{b-1}}\le c \delta^{1-\frac{1}{r}-\e}\prod_{i=1}^3\n{u_i}{\X{r}{s}{b}},\,\,\,\]
provided $2 \ge r > 1,s\ge s(r), b > \frac{1}{r}, \e > 0$. Inserting this estimate, especially the specific power of $\delta$, into the proof of the local result, we 
obtain a lifespan of size $\delta \sim \|u_0\|^{-\frac{2r}{r-1}-\e'}_{\h{r}{s}}$. For $r=2$ this coincides - 
up to $\e'$ - with the result in \cite{KPV93} (see also \cite[Theorem 1.1]{FLP99}).

\item[ii)] Concerning related results for the one-dimensional cubic NLS and DNLS equations we refer to \cite{G05}.
\end{itemize}

{\bf{Acknowledgement:}} The first author, A. G., wishes to thank the Department of Mathematics at the UPV in Bilbao 
for its kind hospitality during his visit.

\section{Bi- and trilinear Airy estimates}

Throughout this section we consider solutions $u(t)=e^{-t\partial_x^3}u_0$, $v(t)=e^{-t\partial_x^3}v_0$
and $w(t)=e^{-t\partial_x^3}w_0$ of the Airy equation with data $u_0$, $v_0$ and $w_0$, respectively. Certain 
bi- and trilinear expressions involving these solutions will be estimated in the spaces $\widehat{L_x^p}(\widehat{L_t^q})$ 
and $\widehat{L^r_{xt}} := \widehat{L_x^r}(\widehat{L_t^r})$, where
\[\n{f}{\widehat{L_x^q}(\widehat{L_t^p})}:= \left( \int \Big{(} \int |\widehat{f}(\xi, \tau)|^{p'} d \tau \Big{)}^{\frac{q'}{p'}}d\xi \right)^{\frac{1}{q'}}, \,\,\,\,\,\frac{1}{q}+\frac{1}{q'}=\frac{1}{p}+\frac{1}{p'}=1.\]
(Below we will always write $p'$, $q'$ etc. to indicate conjugate H\"older exponents, $\widehat{f}$ or $\F f$ denote the Fourier transform of $f$, while 
for the partial Fourier transform in the space variable the symbol $\F_x$ is used.) We begin with the following bilinear estimate, which we state and prove in a slightly more general version
than actually needed.

\begin{lemma}\label{B} 
Let $I^s$ denote the Riesz potential of order $-s$ and let $I^s_-(f,g)$ be defined by its Fourier transform 
(in the space variable):
\[\F_x I_-^s (f,g) (\xi) := \int_* d\xi_1|\xi_1-\xi_2|^s \F_xf(\xi_1)\F_xg(\xi_2),\]
where $\int_*$ is shorthand for $\int_{\xi_1+\xi_2=\xi}$. Then we have
\[\n{I^{\frac{1}{p}}I_-^{\frac{1}{p}}(u,v)}{\widehat{L_x^q}(\widehat{L^p_t})}\le c \n{u_0}{\widehat{L^{r_1}_x}}\n{v_0}{\widehat{L^{r_{2}}_x}},\]
provided $1 \le q \le r_{1,2} \le p \le \infty$ and $\frac{1}{p}+\frac{1}{q}=\frac{1}{r_1}+\frac{1}{r_2}$.
\end{lemma}

Proof: Taking the Fourier transform first in space and then in time we obtain
\[ \F_x I^{\frac{1}{p}}I_-^{\frac{1}{p}}(u,v)(\xi,t)= c|\xi|^{\frac{1}{p}} \int_* d\xi_1|\xi_1 -\xi_2|^{\frac{1}{p}}e^{it(\xi_1^3+\xi_2^3)}\F_xu_0(\xi_1)\F_xv_0(\xi_2)\]
and
\[ \F I^{\frac{1}{p}}I_-^{\frac{1}{p}}(u,v)(\xi,\tau)= c|\xi|^{\frac{1}{p}} \int_* d\xi_1|\xi_1 -\xi_2|^{\frac{1}{p}} \delta(\tau - \xi_1^3- \xi_2^3)\F_xu_0(\xi_1)\F_xv_0(\xi_2),\]
respectively. We use $\delta(g(x)) = \sum_n \frac{1}{|g'(x_n)|} \delta(x-x_n)$, where the sum is taken over all simple zeros of $g$, 
which in our case is
\[g(\xi_1)=\tau - \xi_1 ^3 - \xi_2^3 = \tau - \xi^3 + 3\xi \xi_1 (\xi - \xi_1)\]
with the zeros
\[\xi_1^{\pm}= \frac{\xi \pm y}{2}, \hspace{1,8cm} y:=2\sqrt{\frac{\tau}{3\xi} - \frac{\xi ^2}{12}}\]
and the derivative
\[g'(\xi_1^{\pm})= 3\xi (\xi - 2\xi_1^{\pm}) = \mp 3\xi y.\]
Hence
\begin{eqnarray}\label{201}
 &&\F I^{\frac{1}{p}}I_-^{\frac{1}{p}}(u,v)(\xi,\tau) \\
 &=& c|\xi|^{-\frac{1}{p'}}y^{-\frac{1}{p'}}\left(\F_xu_0(\frac{\xi + y}{2})\F_xv_0(\frac{\xi -y}{2}) + \F_xu_0(\frac{\xi - y}{2})\F_xv_0(\frac{\xi + y}{2})\right). \nonumber
\end{eqnarray}
Using $d\tau = 3 |\xi|ydy$, we see that the $L_{\tau}^{p'}$ - norm of the first contribution equals
\[\left(\int dy |\F_xu_0(\frac{\xi + y}{2})\F_xv_0(\frac{\xi -y}{2})|^{p'}\right)^{\frac{1}{p'}} = c \left(|\F_xu_0|^{p'}*|\F_xv_0|^{p'}(\xi) \right)^{\frac{1}{p'}}.\]
Now Young's inequality is applied to see that
\[\left( \int d\xi (|\F_x{u_0}|^{p'}*|\F_x{v_0}|^{p'}(\xi))^{\frac{q'}{p'}}\right)^{\frac{1}{q'}} \le c \n{u_0}{\widehat{L^{r_1}_x}}\n{v_0}{\widehat{L^{r_{2}}_x}}\]
(cf. the proof of \cite[Lemma 1]{G05}), which is the desired bound. Finally we observe that the second contribution in (\ref{201}) 
can be treated in precisely the same manner with $r_1$ and $r_2$ interchanged.

$\hfill \Box$

Arguing similarly as in the proof of Lemma 2.1 in \cite{G04} we obtain:

\begin{kor}\label{B1} For $p$, $q$, $r_{1,2}$ as in the previous lemma and $b_i > \frac{1}{r_i}$ the estimate
\[\n{I^{\frac{1}{p}}I_-^{\frac{1}{p}}(u_1, u_2)}{\widehat{L^q_x}(\widehat{L^p_t})} \le c \n{u_1}{\X{r_1}{0}{b_1}}\n{u_2}{\X{r_2}{0}{b_2}}\]
is valid.
\end{kor}

The next step is to dualize the preceding corollary. For that purpose we recall the bilinear operator $I^s_+$, 
defined by
\[\F_x I_+^s (f,g) (\xi) := \int_{\xi_1+\xi_2=\xi}d\xi_1|\xi + \xi_2|^s \F_xf(\xi_1)\F_xg(\xi_2),\]
and the linear operators
\[M^s_u v := I_-^s(u,v) \hspace{1cm}\mbox{and}\hspace{1cm}N^s_u w := I_+^s(w,\overline{u}),\]
which are formally adjoint w. r. t. the inner product on $L^2_{xt}$ (cf. \cite[p.3299]{G04}). With this notation,
Corollary \ref{B1} expresses the boundedness of
\[I^{\frac{1}{p}}M^{\frac{1}{p}}_{u_1}: \X{r_2}{0}{b_2} \longrightarrow \widehat{L^q_x}(\widehat{L^p_t})\]
with operator norm $ \le c \n{u_1}{\X{r_1}{0}{b_1}}$. By duality, under the additional hypothesis $1<p,q,r_{1,2}< \infty$, it follows that
\[N^{\frac{1}{p}}_{u_1}I^{\frac{1}{p}}: \widehat{L^{q'}_x}(\widehat{L^{p'}_t}) \longrightarrow \X{r'_2}{0}{-b_2}\]
is bounded with the same norm. Thus we obtain the following estimate:

\begin{kor}\label{B2} Let $1 < q \le r_{1,2} \le p < \infty$, $\frac{1}{p}+\frac{1}{q}=\frac{1}{r_1}+\frac{1}{r_2}$ and 
$b_i > \frac{1}{r_i}$. Then
\begin{equation}\label{202}
\n{I_+^{\frac{1}{p}}(I^{\frac{1}{p}}u_2, \overline{u}_1)}{\X{r'_2}{0}{-b_2}}\le c \n{u_1}{\X{r_1}{0}{b_1}}\n{u_2}{\widehat{L^{q'}_x}(\widehat{L^{p'}_t})}. 
\end{equation}
\end{kor}

\underline{Remark:} Since the phase function $\phi(\xi)= \xi^3$ is odd, we have $\n{u_1}{\x}=\n{\overline{u}_1}{\x}$, 
and we may replace $\overline{u}_1$ by $u_1$ in the left hand side of (\ref{202}).

\vspace{0.5cm}

The special case in (\ref{202}), where $p=q=r_{1,2}$, will be sufficient for our purposes. In this case, (\ref{202}) 
can be written as
\begin{equation}\label{203}
\n{I_+^{\frac{1}{r'}}(I^{\frac{1}{r'}}u_2, u_1)}{\X{{r}}{0}{b'}}\le c \n{u_1}{\X{r'}{0}{-b'}}\n{u_2}{\widehat{L^{r}_{xt}}}, 
\end{equation}
provided $1 < r < \infty$, $b'<-\frac{1}{r'}$. Combining this with the trivial endpoint of the Hausdorff-Young 
inequality, i. e.
\[\n{u_2 u_1}{\widehat{L^{r}_{xt}}}\le c \n{u_1}{\widehat{L^{\infty}_{xt}}}\n{u_2}{\widehat{L^{r}_{xt}}},\]
we obtain by elementary H\"older estimates
\begin{equation}\label{204}
\n{I_+^{\frac{1}{\rho'}}(I^{\frac{1}{\rho'}}u_2, u_1)}{\X{{r}}{0}{\beta}}\le c \n{u_1}{\X{\rho'}{0}{-\beta}}\n{u_2}{\widehat{L^{r}_{xt}}}, 
\end{equation}
where $0 \le \frac{1}{\rho'}\le \frac{1}{r'}$ and $\beta < -\frac{1}{\rho'}$. In this form actually we shall make use of Corollary \ref{B2}.

\vspace{0.5cm}

Now we turn to the trilinear estimates. Again we take the Fourier transform first in $x$ and then in $t$ to 
obtain
\[\F_x (uvw)(\xi, t) = c \int_* d\xi_1 d\xi_2 e^{it(\xi_1^3+\xi_2^3+\xi_3^3)}\F_xu_0(\xi_1)\F_xv_0(\xi_2)\F_xw_0(\xi_3)\]
(where now $\int_*=\int_{\xi_1+\xi_2+\xi_3=\xi}$) and  
\[\F (uvw)(\xi, \tau) = c \int_* d\xi_1 d\xi_2 \delta(\xi_1^3+\xi_2^3+\xi_3^3 - \tau)\F_xu_0(\xi_1)\F_xv_0(\xi_2)\F_xw_0(\xi_3).\]
Now the argument of $\delta$, that is
\[g(\xi_2)=3(\xi-\xi_1)\xi_2^2 - 3(\xi-\xi_1)^2\xi_2 - 3\xi\xi_1(\xi-\xi_1) + \xi^3-\tau\]
has exactly two zeros
\begin{equation}\label{205}
\xi_2^{\pm}=\frac{\xi-\xi_1}{2}\pm \sqrt{\frac{(\xi+\xi_1)^2}{4}+\frac{\tau-\xi^3}{3(\xi-\xi_1)}}=:\frac{\xi-\xi_1}{2}\pm y,
\end{equation}
with
\[|g'(\xi_2^{\pm})|=6|\xi-\xi_1|\sqrt{\frac{(\xi+\xi_1)^2}{4}+\frac{\tau-\xi^3}{3(\xi-\xi_1)}}=6|\xi-\xi_1|y.\]
Using $\delta(g(\xi_2)) = \sum_{g(x_n)=0} \frac{\delta(\xi_2 - x_n)}{|g'(x_n)|}$, where the sum is taken over all simple zeros of $g$,
we see that
\[\F (uvw)(\xi, \tau) = c(K_+(\xi, \tau)+K_-(\xi, \tau)),\]
where
\[K_{\pm}(\xi, \tau)= \int d\xi_1 \frac{1}{|\xi-\xi_1|y}\F_xu_0(\xi_1)\F_xv_0(\frac{\xi-\xi_1}{2}\pm y)\F_xw_0(\frac{\xi-\xi_1}{2}\mp y)\]
with $y$ as defined in (\ref{205}).

In order to estimate $\n{uvw}{\widehat{L^r_{xt}}}=\n{\F (uvw)}{L^{r'}_{\xi \tau}}$ we distinguish between three 
cases depending on the relative size of the frequencies $\xi_1, \xi_2$ and $\xi_3$:

\begin{itemize}
\item[i)] $|\xi_1| \sim |\xi_2| \gg \langle\xi_3\rangle$,
\item[ii)] $|\xi_2-\xi_3|\ge|\xi_2+\xi_3|$,
\item[iii)] $1\le|\xi_2-\xi_3|\le|\xi_2+\xi_3|$.
\end{itemize}

To treat the first case we define the trilinear operator $T$ by
\[\F_x T(f,g,h):= \int_* d\xi_1  d\xi_2 \F_xf(\xi_1)\F_xg(\xi_2)\F_xh(\xi_3)\chi_{\{|\xi_1| \sim |\xi_2| \gg \langle\xi_3\rangle\}},\]
where again $\int_*=\int_{\xi_1+\xi_2+\xi_3=\xi}$. In this case we have:

\begin{lemma}\label{T1}
Let $1\le r\le2$ and $s_1>\frac{1}{4r'}-\frac{1}{2}$, $s_2 \ge\frac{1}{2r'}$. Then
\[\n{T(u,v,w)}{\widehat{L^r_{xt}}} \le c \n{u_0}{\widehat{H^r_{s_1}}}\n{v_0}{\widehat{H^r_{s_1}}}\n{w_0}{\widehat{H^r_{s_2}}}.\]
\end{lemma}

Proof: By the above computation we have
\[\F T(u,v,w)(\xi,\tau)= c (K^+(\xi,\tau) + K^-(\xi,\tau)),\]
with
\[K^{\pm}(\xi,\tau)=\int_{A_{\pm}} d\xi_1 \frac{1}{|\xi-\xi_1|y}\F_xu_0(\xi_1)\F_xv_0(\frac{\xi-\xi_1}{2}\pm y)\F_xw_0(\frac{\xi-\xi_1}{2}\mp y),\]
where $A_{\pm}=\{|\xi_1| \sim |\frac{\xi-\xi_1}{2}\pm y| \gg \langle\frac{\xi-\xi_1}{2}\mp y\rangle\}$ and $y$ is defined by (\ref{205}). 
Since in $A_{\pm}$ the inequality $|\xi_1||\frac{\xi-\xi_1}{2}\pm y| \le c |\xi-\xi_1|y$ holds true, we get the upper bound 
\[K^{\pm}(\xi,\tau) \le c \int d\xi_1 \F_x J^{-1}u_0(\xi_1)\F_xJ^{-1}v_0(\frac{\xi-\xi_1}{2}\pm y)\F_xw_0(\frac{\xi-\xi_1}{2}\mp y),\]
leading to
\[\n{T(u,v,w)}{\widehat{L^1_{xt}}} \le c \n{J^{-1}u_0}{\widehat{L^{\infty}_{x}}}\n{J^{-1}v_0}{\widehat{L^1_{x}}}\n{w_0}{\widehat{L^1_{x}}}.\]
By symmetry between the first two factors and multilinear interpolation we obtain
\begin{equation}\label{206}
\n{T(u,v,w)}{\widehat{L^1_{xt}}} \le c \n{J^{-1}u_0}{L^{2}_{x}}\n{J^{-1}v_0}{L^2_{x}}\n{w_0}{\widehat{L^1_{x}}}.
\end{equation}
On the other hand side we have
\[\n{uvw}{L^2_{xt}} \le c \n{u}{L_x^8(L_t^4)}\n{v}{L_x^8(L_t^4)}\n{w}{L_x^4(L_t^{\infty})},\]
with
\begin{equation}\label{207}
\n{w}{L_x^4(L_t^{\infty})} \le c \n{I^{\frac{1}{4}}u_0}{L_x^2},
\end{equation}
which is the maximal function estimate from \cite[Thm. 3]{S87}. Concerning the first two factors we interpolate between 
the sharp version of Kato's smoothing effect, \newline i. e. $\n{Iu}{L_x^{\infty}(L_t^2)} = c \n{u_0}{L_x^2}$, see \cite[Thm. 4.1]{KPV91},
and (\ref{207}) to obtain
\[\n{I^{\frac{3}{8}}u}{L_x^8(L_t^4)}\le c \n{u_0}{L_x^2},\]
such that
\begin{equation}\label{208}
\n{T(u,v,w)}{L^2_{xt}} \le c \n{J^{-\frac{3}{8}}u_0}{L_x^2}\n{J^{-\frac{3}{8}}v_0}{L_x^2}\n{J^{\frac{1}{4}}w_0}{L_x^2}.
\end{equation}
Using multilinear interpolation again, now between (\ref{206}) and (\ref{208}), we finally see that, for $1\le r\le 2$,
\begin{eqnarray*}
\n{T(u,v,w)}{\widehat{L^r_{xt}}} & \le & c\n{J^{\frac{5}{4r'}-1}u_0}{L_x^2}\n{J^{\frac{5}{4r'}-1}v_0}{L_x^2}\n{J^{\frac{1}{2r'}}w_0}{\widehat{L^r_{x}}}\\
& \le & c \n{u_0}{\widehat{H^r_{s_1}}}\n{v_0}{\widehat{H^r_{s_1}}}\n{w_0}{\widehat{H^r_{s_2}}},
\end{eqnarray*}
where in the last step we have used the Sobolev type embedding $\h{r}{s} \subset \h{\rho}{\sigma}$, which holds true for $s-\frac{1}{r}>\sigma-\frac{1}{\rho}$, $r\le \rho$.

$\hfill \Box$

\begin{kor}\label{T1C}
For $r$, $s_{1,2}$ as in the previous lemma and $b>\frac{1}{r}$ the estimate
\[\n{T(u_1,u_2,u_3)}{\widehat{L^r_{xt}}} \le c \n{u_1}{\X{r}{s_1}{b}}\n{u_2}{\X{r}{s_1}{b}}\n{u_3}{\X{r}{s_2}{b}}\]
holds true.
\end{kor}

Next we introduce $T_{\ge}$ ($T_{\le}$) by
\[\F_x T_{\ge}(f,g,h):= \int_* d\xi_1  d\xi_2 \F_xf(\xi_1)\F_xg(\xi_2)\F_xh(\xi_3)\chi_{\{|\xi_2-\xi_3| \ge |\xi_2+\xi_3|\}},\]
and
\[\F_x T_{\le}(f,g,h):= \int_* d\xi_1  d\xi_2 \F_xf(\xi_1)\F_xg(\xi_2)\F_xh(\xi_3)\chi_{\{1\le|\xi_2-\xi_3| \le |\xi_2+\xi_3|\}}.\]

\begin{lemma}\label{T2}
Let $1<p_1<p<p_0<\infty$, $p<p'_0$, $\frac{3}{p}=\frac{1}{p_0}+\frac{2}{p_1}$ and $\frac{2}{p_1}<1+\frac{1}{p}$. Then the estimate
\[\n{T_{\ge}(u,v,w)}{\widehat{L^p_{xt}}} \le c \n{u_0}{\widehat{L^{p_0}_{x}}}\n{I^{-\frac{1}{2p}}v_0}{\widehat{L^{p_1}_{x}}}\n{I^{-\frac{1}{2p}}w_0}{\widehat{L^{p_1}_{x}}}\]
is valid.
\end{lemma}

Proof: For the Fourier transform of $T_{\ge}(u,v,w)$ in both variables we obtain
\[\F T_{\ge}(u,v,w)(\xi,\tau)= c (K^+_{\ge}(\xi,\tau) + K^-_{\ge}(\xi,\tau)),\]
where
\[K^{\pm}_{\ge}(\xi,\tau)=\int_{\{2y\ge|\xi-\xi_1|\}} d\xi_1 \frac{1}{|\xi-\xi_1|y}\F_xu_0(\xi_1)\F_xv_0(\frac{\xi-\xi_1}{2}\pm y)\F_xw_0(\frac{\xi-\xi_1}{2}\mp y),\]
with $y$ as in (\ref{205}) again. By symmetry we may restrict ourselves to the estimation of $K^+_{\ge}$. Using $|\frac{\xi-\xi_1}{2}\pm y|\le 2y$ and H\"older's inequality, we
see that
\begin{eqnarray*}
&&K^+_{\ge}(\xi,\tau) \le c \left( \int d \xi_1 \frac{|\F_x u_0 (\xi_1)|^p}{|\xi-\xi_1|^{(1-\theta)p}}\right)^{\frac{1}{p}}\times \dots \\
\dots & \times & \left( \int \frac{d \xi_1}{|\xi-\xi_1|^{\theta p'}y}|\F_xI^{-\frac{1}{2p}}v_0(\frac{\xi-\xi_1}{2}+ y)\F_xI^{-\frac{1}{2p}}w_0(\frac{\xi-\xi_1}{2}- y)|^{p'}\right)^{\frac{1}{p'}},
\end{eqnarray*}
where $\theta = \frac{3}{p'}-\frac{2}{p_1'}$ ($\in (0,1)$ by our assumptions). Taking the $L^{p'}_{\tau}$-norm 
of both sides and using $d\tau = 6 |\xi-\xi_1|ydy$ we arrive at
\begin{eqnarray*}
&& \n{\F T_{\ge}(u,v,w)(\xi,\cdot)}{L^{p'}_{\tau}} \le c (|\F_x u_0 |^p * |\xi|^{(\theta - 1)p})^{\frac{1}{p}}\times \dots \\
\dots & \times & \left( \int \frac{d \xi_1 dy}{|\xi-\xi_1|^{\theta p' - 1}}|\F_xI^{-\frac{1}{2p}}v_0(\frac{\xi-\xi_1}{2}+ y)\F_xI^{-\frac{1}{2p}}w_0(\frac{\xi-\xi_1}{2}- y)|^{p'}\right)^{\frac{1}{p'}}.
\end{eqnarray*}
Changing variables ($z_{\pm}:=\frac{\xi-\xi_1}{2}\pm y$) we see that the second factor equals
\[\left(\!\int\! \frac{d z_+ dz_-}{|z_++z_-|^{\theta p' - 1}}|\F_xI^{-\frac{1}{2p}}v_0(z_+)\F_xI^{-\frac{1}{2p}}w_0(z_-)|^{p'}\!\right)^{\frac{1}{p'}} \!\!\!\le c\n{I^{-\frac{1}{2p}}v_0}{\widehat{L^{p_1}_{x}}} \n{I^{-\frac{1}{2p}}w_0}{\widehat{L^{p_1}_{x}}},\]
by the Hardy-Littlewood-Sobolev-inequality, requiring $\theta$ to be chosen as above and $1<\theta p'<2$, which follows from our assumptions.
It remains to estimate the $L^{p'}_{\xi}$-norm of the first factor, that is
\begin{eqnarray*}
&& \||\F_x{u_0}|^p* |\xi|^{(\theta -1)p}\|^{\frac{1}{p}}_{L_{\xi}^{\frac{p'}{p}}}\\
& \le & c( \n{|\F_x{u_0}|^p}{L_{\xi}^{\frac{p'_0}{p}}} \n{|\xi|^{(\theta -1)p}}{L_{\xi}^{\frac{1}{(1-\theta)p}, \infty}} )^{\frac{1}{p}}\\
& \le & c \n{u_0}{\widehat{L^{p_0}_x}},
\end{eqnarray*}
where the HLS inequality was used again. For its application we need
\[0<(1-\theta)p<1;\,\,\,\,1<\frac{p'_0}{p}< \frac{1}{1-(1-\theta)p}\,\,\,\,\mbox{and}\,\,\,\,\,\theta=\frac{1}{p'_0},\]
which follows from the assumptions, too.
$\hfill \Box$

\begin{kor}\label{T2C}
For $1<r<2$ there exist $s_{0,1}\ge 0$ with $s_0 + 2s_1 = \frac{1}{r}$, such that
\begin{equation}\label{209}
\n{T_{\ge}(u,v,w)}{\widehat{L^r_{xt}}} \le c \n{I^{-s_0}u_0}{\widehat{L^{r}_{x}}}\n{I^{-s_1}v_0}{\widehat{L^{r}_{x}}}\n{I^{-s_1}w_0}{\widehat{L^{r}_{x}}}.
\end{equation}
In addition, for $b>\frac{1}{r}$ we have
\[\n{T_{\ge}(u_1,u_2,u_3)}{\widehat{L^r_{xt}}} \le c \n{I^{-s_0}u_1}{\X{r}{0}{b}}\n{I^{-s_1}u_2}{\X{r}{0}{b}}\n{I^{-s_1}u_3}{\X{r}{0}{b}}.\]
\end{kor}

Proof of (\ref{209}): Using H\"older's inequality and the Airy-version of the Fefferman-Stein-estimate, that is
\begin{equation}\label{FS}
\n{u}{L^{3q}_{xt}} \le c \n{I^{-\frac{1}{3q}}u_0}{\widehat{L^q_{x}}},\hspace{1cm}q>\frac{4}{3},
\end{equation}
see \cite[Corollary 3.6]{G04}, we get for
\begin{equation}\label{210}
\frac{4}{3}< q_0<2<q_1,\hspace{1cm}\mbox{with}\hspace{1cm}\frac{3}{2}=\frac{1}{q_0}+\frac{2}{q_1}
\end{equation}
that
\begin{equation}\label{211}
\n{T_{\ge}(u,v,w)}{L^{2}_{xt}} \le \n{uvw}{L^{2}_{xt}} \le c \n{I^{-\frac{1}{3q_0}}u_0}{\widehat{L^{q_0}_{x}}}\n{I^{-\frac{1}{3q_1}}v_0}{\widehat{L^{q_1}_{x}}}\n{I^{-\frac{1}{3q_1}}w_0}{\widehat{L^{q_1}_{x}}}.
\end{equation}
Multilinear interpolation of (\ref{211}) with Lemma \ref{T2} yields (\ref{209}), provided $p,p_0,p_1$; $q_0,q_1$, defined by the interpolation conditions
\[\frac{1}{r}=\frac{1- \theta}{p}+\frac{\theta}{2}=\frac{1-\theta}{p_0}+\frac{\theta}{q_0}=\frac{1-\theta}{p_1}+\frac{\theta}{q_1},\]
fulfill the assumptions of Lemma \ref{T2} and (\ref{210}), respectively, which can be guaranteed by choosing $\theta$ 
sufficiently small. Now $s_{0,1}$ are obtained from
\[s_0=\frac{\theta}{3q_0}\hspace{1cm}\mbox{and}\hspace{1cm}s_1=\frac{1-\theta}{2p}+\frac{\theta}{3q_1},\]
which gives
\[s_0 + 2s_1 =\frac{1-\theta}{p}+ \frac{\theta}{3}(\frac{1}{q_0}+\frac{2}{q_1}) = \frac{1}{r}\]
as desired.

$\hfill \Box$

\vspace{0.1cm}

\underline{Remark:} By (\ref{FS}), Corollary \ref{T2C} still holds true for $r \ge 2$ (with $s_0=s_1=\frac{1}{3r}$).

\vspace{0.1cm}

\begin{lemma}\label{T3}
Let $1 \le r < \rho \le \infty$. Then
\[\n{T_{\le}(u,v,w)}{\widehat{L^r_{xt}}} \le c \n{u_0}{\widehat{L^{\rho}_{x}}}\n{I^{-\frac{1}{2r}}v_0}{\widehat{L^{r}_{x}}}\n{I^{-\frac{1}{2r}}w_0}{\widehat{L^{r}_{x}}}.\]
\end{lemma}

Proof: We have
\[\F T_{\le}(u,v,w)(\xi,\tau)= c (K^+_{\le}(\xi,\tau) + K^-_{\le}(\xi,\tau)),\]
where
\[K^{\pm}_{\le}(\xi,\tau)=\int_{\{1 \le 2y \le|\xi-\xi_1|\}}  \frac{d\xi_1}{|\xi-\xi_1|y}\F_xu_0(\xi_1)\F_xv_0(\frac{\xi-\xi_1}{2}\pm y)\F_xw_0(\frac{\xi-\xi_1}{2}\mp y)\]
with $y$ as defined in (\ref{205}).
By symmetry between $v$ and $w$ it suffices to treat $K^+_{\le}$, which we decompose dyadically with respect to $y$ to obtain the upper bound:
\begin{eqnarray*}
\!\!\!\!\!&  & c \sum_{j=0}^{\infty}\int_{\{1 \le 2y \le|\xi-\xi_1|\,, \,\,y \sim 2^j\}}  \frac{d\xi_1}{|\xi-\xi_1|y}\F_xu_0(\xi_1)\F_xv_0(\frac{\xi-\xi_1}{2}+ y)\F_xw_0(\frac{\xi-\xi_1}{2}- y) \\
\!\!\!\!\!& \le & c \sum_{j=0}^{\infty}2^{-j}\int_{\{y \sim 2^j\}}d\xi_1\F_xu_0(\xi_1)\F_xI^{-\frac{1}{2}}v_0(\frac{\xi-\xi_1}{2}+ y)\F_xI^{-\frac{1}{2}}w_0(\frac{\xi-\xi_1}{2}- y) \\
\!\!\!\!\!& \le & c \sum_{j=0}^{\infty}2^{-j}\n{u_0}{\widehat{L^{p}_{x}}}\lambda(\{y \sim 2^j\})^{\frac{1}{p'}}\n{I^{-\frac{1}{2}}v_0}{\widehat{L^{1}_{x}}}\n{I^{-\frac{1}{2}}w_0}{\widehat{L^{1}_{x}}},
\end{eqnarray*}
where $\lambda(\{y \sim 2^j\})$ denotes the Lebesgue measure of $\{\xi_1 : y(\xi_1) \sim 2^j\}$, which is bounded by $c2^j$\footnote{To see this,
we write $\{\xi_1 : y(\xi_1) \sim 2^j\}= S_1 \cup S_2$, where in $S_1$ we assume that $|\xi-\xi_1| \lesssim 2^j $, 
$|\xi+\xi_1| \lesssim 2^j $ or $|\xi-3\xi_1| \lesssim 2^j $. Then $S_1$ consists of a finite number of intervals of total length 
bounded by $c2^j$. For $S_2$ we have $|\xi-\xi_1| \gg 2^j$, $|\xi+\xi_1| \gg 2^j$ and $|\xi-3\xi_1| \gg 2^j$, implying that
\[\left|\frac{dy}{d\xi_1}\right| = \frac{1}{2y|\xi-\xi_1|}\left|\frac{(\xi+\xi_1)(\xi-3\xi_1)}{4}+ y^2\right| \gtrsim \frac{|\xi+\xi_1||\xi-3\xi_1|}{y|\xi-\xi_1|} \gtrsim 1\,\,\,,\]
which gives
\[\lambda (S_2) = \int_{S_2} d\xi_1  \le \int \frac{d\xi_1}{dy}\chi_{\{y \sim 2^j\}} dy \le c 2^j.\] }.
Hence, for any $p>1$,
\begin{eqnarray}\label{220}
\n{K^+_{\le}}{L^{\infty}_{\xi \tau  }} & \le & c \sum_{j=0}^{\infty}2^{-\frac{j}{p'}}\n{u_0}{\widehat{L^{p}_{x}}}\n{I^{-\frac{1}{2}}v_0}{\widehat{L^{1}_{x}}}\n{I^{-\frac{1}{2}}w_0}{\widehat{L^{1}_{x}}} \nonumber \\
& \le & c \n{u_0}{\widehat{L^{p}_{x}}}\n{I^{-\frac{1}{2}}v_0}{\widehat{L^{1}_{x}}}\n{I^{-\frac{1}{2}}w_0}{\widehat{L^{1}_{x}}}.
\end{eqnarray}
On the other hand, by integration with respect first to $d\tau = 6y(\xi-\xi_1)dy$, to $d\xi$ and finally to $d\xi_1$, 
we see that
\begin{equation}\label{221}
\n{K^+_{\le}}{L^{1}_{\xi \tau  }} \le c \n{u_0}{\widehat{L^{\infty}_{x}}}\n{v_0}{\widehat{L^{\infty}_{x}}}\n{w_0}{\widehat{L^{\infty}_{x}}}.
\end{equation}
Now multilinear interpolation between (\ref{220}) and (\ref{221}) leads to
\[\n{K^+_{\le}}{L^{r'}_{\xi \tau  }} \le c \n{u_0}{\widehat{L^{\rho}_{x}}}\n{I^{-\frac{1}{2r}}v_0}{\widehat{L^{r}_{x}}}\n{I^{-\frac{1}{2r}}w_0}{\widehat{L^{r}_{x}}},\]
which gives the desired result.
$\hfill \Box$

\begin{kor}\label{T3C}
Let $1\le r < \rho \le \infty$, $\beta>\frac{1}{\rho}$, $b>\frac{1}{r}$ and $\e>0$. Then
\[\n{T_{\le}(u_1,u_2,u_3)}{\widehat{L^r_{xt}}} \le c \n{u_1}{\X{\rho}{0}{\beta}}\n{I^{-\frac{1}{2r}}u_2}{\X{r}{0}{b}}\n{I^{-\frac{1}{2r}}u_3}{\X{r}{0}{b}}\]
and
\[\n{T_{\le}(u_1,u_2,u_3)}{\widehat{L^r_{xt}}} \le c \n{u_1}{\X{r}{\e}{b}}\n{I^{-\frac{1}{2r}}u_2}{\X{r}{0}{b}}\n{I^{-\frac{1}{2r}}u_3}{\X{r}{0}{b}}\]
are valid.
\end{kor}

\section{Proof of Theorem 2}

Without loss of generality we may assume that $s = s(r)$. Then we rewrite the left hand side of (\ref{40}) as
\[\n{\langle \tau -  \xi^3 \rangle^{b'} \langle  \xi \rangle^{s}|\xi|\int d \nu  \prod_{i=1}^3 \widehat{u_i}(\xi_i,\tau_i)}{L^{r'}_{\xi,\tau}},\]
where $d \nu = d\xi_1 d\xi_{2} d\tau_1 d \tau_{2}$ and $\sum_{i=1}^3 (\xi_i,\tau_i) = (\xi, \tau)$. \\

In the sequel, we shall use the following notation:

\begin{itemize}
\item $\xi_{max}$, $\xi_{med}$, $\xi_{min}$ are defined by $|\xi_{max}| \ge |\xi_{med}| \ge |\xi_{min}|$,
\item $p$ denotes the projection on low frequencies, i. e. $\widehat{pf}(\xi)=\chi_{\{|\xi| \le 1\}}\widehat{f}(\xi)$,
\item $f \preceq g$ is shorthand for $|\widehat{f}| \le c |\widehat{g}|$,
\item for the mixed weights coming from the $\x$ - norms we shall write $\sigma_0 := \tau - \xi^3$ and $\sigma_i := \tau_i - \xi_i^3$, $1 \le i \le 3$, respectively,
\item the Fourier multiplier associated with these weights is denoted by $\Lambda ^b := \F ^{-1} \langle \tau - \xi^3 \rangle ^b \F$,
\item for a real number $x$ we write $x\pm$ to denote $x \pm \e$ for arbitrarily small $\e > 0$, $\infty -$ stands for an arbitrarily large real number.
\end{itemize}

Apart from the trivial region where $|\xi_{max}| \le 1$, whose contribution can be estimated by
\[\n{\prod_{i=1}^3 p u_i}{\widehat{L^r_{xt}}}\le c \prod_{i=1}^3 \n{p u_i}{\widehat{L^{3r}_{xt}}}\le c \prod_{i=1}^3 \n{p u_i}{\X{r}{0}{b}}\le c \prod_{i=1}^3 \n{u_i}{\x},\]
we consider three cases:

\begin{itemize}
\item[1.] The nonresonant case, where $|\xi_{max}| \gg |\xi_{med}|$,
\item[2.] the semiresonant case with $|\xi_{max}| \sim |\xi_{med}| \gg |\xi_{min}|$
and, finally,
\item[3.] the resonant case, where $|\xi_{max}| \sim |\xi_{min}|$.
\end{itemize}

\vspace{0.5cm}

1. In the nonresonant case we assume without loss of generality that $|\xi_1|\ge|\xi_2|\ge|\xi_3|$. Then 
we have for this region
\begin{eqnarray*}
J^s \partial_x(u_1u_2u_3) & \preceq & \partial_x(J^su_1J^su_2J^{-s}u_3) \\
& \preceq & I^{\frac{1}{r}}I_-^{\frac{1}{r}}(J^su_1,J^su_2)J^{1-s-\frac{2}{r}}u_3 \\
& \preceq & I_+^{0+}(I^{\frac{1}{r}+}I_-^{\frac{1}{r}}(J^su_1,J^su_2),J^{1-s-\frac{2}{r}-}u_3) .
\end{eqnarray*}
Now the dual version (\ref{204}) of the bilinear estimate is applied to obtain
\begin{eqnarray*}
&& \n{I_+^{0+}(I^{\frac{1}{r}+}I_-^{\frac{1}{r}}(J^su_1,J^su_2),J^{1-s-\frac{2}{r}-}u_3)}{\X{r}{0}{b'}}\\
& \le & c \n{I^{\frac{1}{r}}I_-^{\frac{1}{r}}(J^su_1,J^su_2)}{\widehat{L^r_{xt}}} \n{J^{1-s-\frac{2}{r}-}u_3}{\X{\infty-}{0}{0+}}\le c \prod_{i=1}^3 \n{u_i}{\x},
\end{eqnarray*}
where in the last step we have used the bilinear estimate itself (Corollary \ref{B1}) for the first and Sobolev-type embeddings for the second factor.

\vspace{0.5cm}

2. In the semiresonant case we assume again $|\xi_1|\ge|\xi_2|\ge|\xi_3|$ and consider two subcases: If,
in addition, $|\xi_1+\xi_2| \le 1$ (so that $\langle{\xi}\rangle \le c \langle{\xi_3}\rangle$), we can argue 
as in case 1, with $u_1$ and $u_3$ interchanged:
\begin{eqnarray*}
J^s \partial_x(u_1u_2u_3) & \preceq & \partial_x(J^{-s}u_1J^su_2J^{s}u_3) \\
\preceq \dots & \preceq & I_+^{0+}(I^{\frac{1}{r}+}I_-^{\frac{1}{r}}(J^su_3,J^su_2),J^{1-s-\frac{2}{r}-}u_1),
\end{eqnarray*}
which can be treated as above by applying (\ref{204}), Sobolev-type embeddings and Corollary \ref{B1}. On the other hand, if $|\xi_1+\xi_2| \ge 1$, we have
\[|\sigma_0 - \sigma_1 - \sigma_2 - \sigma_3| = 3|\xi_1+\xi_2||\xi_2+\xi_3||\xi_3+\xi_1| \gtrsim \langle \xi_1 \rangle \langle \xi_2 \rangle,\]
and hence, for any $\e>0$,
\[\langle \xi_1 \rangle^{\e} \langle \xi_2 \rangle^{\e} \le c \prod_{i=0}^3 \langle \sigma_i \rangle^{\e}.\]
So, in this subcase, we have the upper bound
\[\n{T(J^{\frac{s+1}{2}-}\Lambda^{0+}u_1,J^{\frac{s+1}{2}-}\Lambda^{0+}u_2,\Lambda^{0+}u_3)}{\widehat{L^r_{xt}}}   \le c \prod_{i=1}^3 \n{u_i}{\x}\]
by Corollary \ref{T1C}.

\vspace{0.5cm}

3. In the resonant case we distinguish several subcases:

\vspace{0.3cm}
\underline{3.1:} At least for one pair $(i,j)$ we have $|\xi_i-\xi_j|\ge|\xi_i+\xi_j|$. \\

Here we may assume by symmetry that $|\xi_2-\xi_3|\ge|\xi_2+\xi_3|$. Then we have for nonnegative $s_{0,1}$ with $s_0 + 2s_1 = \frac{1}{r}$
\[\partial_x J^s (u_1u_2u_3) \preceq T_{\ge}(J^{s+s_0}u_1,J^{s+s_1}u_2,J^{s+s_1}u_3),\]
so that Corollary \ref{T2C} leads to the desired bound.  

\vspace{0.3cm}
\underline{3.2:} $|\xi_1-\xi_2|\le|\xi_1+\xi_2|$, $|\xi_2-\xi_3|\le|\xi_2+\xi_3|$ and $|\xi_3-\xi_1|\le|\xi_3+\xi_1|$, so that all the 
$\xi_i$ have the same sign, which implies
\[|\xi_1|^3 \sim |\xi_2|^3 \sim |\xi_3|^3 \le \prod_{i=0}^3 \langle \sigma_i \rangle.\] 

\underline{3.2.1:} At least one of the $|\xi_i-\xi_j|\ge 1$. \\

By symmetry we may assume that $|\xi_2-\xi_3|\ge 1$. Gaining a $\langle \xi \rangle ^{\e}$ from the $\sigma's$ we obtain as an upper bound for this subcase
\[ \n{T_{\le}(J^{s-}\Lambda^{0+}u_1,J^{\frac{1}{2}}\Lambda^{0+}u_2,J^{\frac{1}{2}}\Lambda^{0+}u_3)}{\widehat{L^r_{xt}}}   \le c \prod_{i=1}^3 \n{u_i}{\x},\]
where we have used the second part of Corollary \ref{T3C}. \\

\underline{3.2.2:} $|\xi_i-\xi_j|\le 1$ for all $1 \le i \neq j \le 3$. \\

Again, we can gain a $\langle \xi \rangle ^{\e}$ from the $\sigma's$. Now, writing 
\[f_i(\xi,\tau)=\langle \xi \rangle^s \langle \tau- \xi^3\rangle^b \F u_i(\xi,\tau), \,\,\,\,1\le i\le 3, \,\,\,\,\mbox{such that}\,\,\,\,
\n{f_i}{L^{r'}_{\xi \tau}}=\n{u_i}{\x},\] 
it suffices to show 
\begin{equation}\label{300}
\n{\langle \xi \rangle^{s-}|\xi| \int_A d\nu \prod_{i=1}^3\langle \xi_i \rangle^{-s}\langle \tau_i- \xi_i^3\rangle^{-\frac{1}{r}-}f_i(\xi_i,\tau_i)}{L^{r'}_{\xi \tau}} \le c \prod_{i=1}^3 \n{f_i}{L^{r'}_{\xi \tau}},
\end{equation}
where in $A$ all the differences $|\xi_k-\xi_j|$, $1 \le k \neq j \le 3$, are bounded by $1$ and $|\xi| \sim |\xi_i| \sim \langle \xi_i \rangle$ for all $1\le i\le 3$.
By H\"older's inequality and Fubini's Theorem the proof of (\ref{300}) is reduced to show that
\begin{equation}\label{301}
\sup_{\xi, \tau}\langle \xi \rangle^{1-2s-} \left( \int_A d\nu \prod_{i=1}^3 \langle \tau_i- \xi_i^3\rangle^{-1-}\right)^{\frac{1}{r}} < \infty .
\end{equation}
Using \cite[Lemma 4.2]{GTV97} twice, we see that
\[\int_A d\nu \prod_{i=1}^3 \langle \tau_i- \xi_i^3\rangle^{-1-} \le c \int_{A'} d\xi_1d\xi_2 \langle \tau- \xi^3 + 3(\xi_1+\xi_2)(\xi-\xi_1)(\xi-\xi_2)\rangle^{-1-},\]
where $A'$ is simply the projection of $A$ onto $\mathbb{R}^2$. We decompose
\[A'=A_0 \cup A_1 \cup \bigcup_{0\le k,j \le c \ln{(|\xi|)}}A_{kj},\]
where in $A_0$ ($A_1$) we have that $|\xi_1+\xi_2-\frac{2\xi}{3}| \le \frac{100}{|\xi|}$ ($|\xi_1+\xi_3-\frac{2\xi}{3}| \le \frac{100}{|\xi|}$), so that the contributions of 
these subregions are bounded by $\frac{c}{|\xi|}$, while in $A_{kj}$ it should hold that $|\xi_1+\xi_2-\frac{2\xi}{3}|\sim 2^{-k}$ and $|\xi_1+\xi_3-\frac{2\xi}{3}|\sim 2^{-j}$.
By symmetry we may assume $k \le j$. To estimate the integral over $A_{kj}$, we introduce new variables $x_1:=\xi_1+\xi_2-\frac{2\xi}{3}$ and $x_2:=\xi_1-\xi_2$, such that
\[|x_1| \sim 2^{-k}\hspace{1cm}\mbox{and}\hspace{1cm}|x_2|=|\xi_1+\xi_2-\frac{2\xi}{3} +2(\xi_1+\xi_3-\frac{2\xi}{3})| \lesssim 2^{-k}.\]
Then
\begin{eqnarray*}
&& \int_{A_{kj}} d\xi_1d\xi_2 \langle \tau- \xi^3 + 3(\xi_1+\xi_2)(\xi-\xi_1)(\xi-\xi_2)\rangle^{-1-} \\
& \le & \int_{|x_2|\lesssim 2^{-k}} \!\!\!\!\!\!\!\!\!\!dx_2 \int_{|x_1|\sim 2^{-k}} \!\!\!\!\!\!\!\!\!\!dx_1 \langle \textstyle{\tau- \xi^3 + 3(x_1+\frac{2\xi}{3})(\frac{x_1+x_2}{2}-\frac{2\xi}{3})(\frac{x_1-x_2}{2}-\frac{2\xi}{3})}\rangle^{-1-}.
\end{eqnarray*}
Substituting $z:=(x_1+\frac{2\xi}{3})(\frac{x_1+x_2}{2}-\frac{2\xi}{3})(\frac{x_1-x_2}{2}-\frac{2\xi}{3})$, so that
\[\left| \frac{dz}{dx_1}\right| = \left| \frac{3x_1^2 - x_2^2}{4} - x_1\xi\right| \sim |x_1\xi|\sim |\xi|2^{-k},\]
we see that the latter is bounded by
\[\int_{|x_2|\lesssim 2^{-k}}dx_2 \int 2^k \frac{dz}{|\xi|}\langle \tau- \xi^3 + 3z\rangle^{-1-} \le \frac{c}{|\xi|}.\]
Finally, summing up over $j$ and $k$, we have
\[\int_{A'} d\xi_1d\xi_2 \langle \tau- \xi^3 + 3(\xi_1+\xi_2)(\xi-\xi_1)(\xi-\xi_2)\rangle^{-1-}\le c \frac{(\ln{|\xi|})^2}{|\xi|} \le c |\xi|^{-1+},\]
which gives (\ref{301}).

$\hfill \Box$

\end{document}